\documentclass[12pt]{article}

\usepackage{amsmath}
\usepackage{amssymb}
\usepackage{amsthm}

\usepackage{array}

\usepackage[margin=1in]{geometry}  %
\usepackage{pdflscape}                %
\usepackage{graphicx}             %
\usepackage{siunitx}
\usepackage{comment}

\usepackage{adjustbox}

\newtheorem{theorem}{Theorem}

\usepackage[backend=biber, style=authoryear]{biblatex}
\newcolumntype{L}{>{$}l<{$}} %
\newcolumntype{C}{>{$}c<{$}} %
\newcolumntype{R}{>{$}r<{$}} %

\newmuskip\pFqskip
\mathchardef\pFcomma=\mathcode`, %

\begin{document}

\title{Intractable enumeration problems are like Russian nesting dolls: structural properties of monomer-dimer coverings on two-dimensional quadratic lattices}

\author{Yong Kong\\
  Yale University \\
  New Haven, CT 06520, USA \\
  email: \texttt{yong.kong@yale.edu}
}%

\date{}
\maketitle

\begin{abstract}
  Counting the number of coverings of $s$ dimers on two-dimensional quadratic lattices is considered as intractable and belongs to \#P-complete class.
  We reveal the structure of the exact solution to the problem and provide an explicit formula for it,
  which includes $s-1$ nesting sums. 
  This results in an exponential time complexity of $O(2^s)$.
  The solution is explicitly determined by a sequence that exhibits double-exponential growth.
\end{abstract}

\section{Introduction}

Physicists have studied rigid hard-core molecules on lattices for a long time
as models for phase transitions in equilibrium statistical mechanics.
For the monomer-polymer models, a linear rigid polymer covers $k$ adjacent lattice sites,
with no lattice site occupied by more than one polymer.
The polymers are usually called $k$-mers,
and those unoccupied lattice sites are called monomers.
As a special case, when $k=2$, the model becomes the well-known monomer-dimer model,
which is closely related to the Ising model.
While
the special case of close packed dimer problem can be solved for any \emph{planar} lattices
using the Pfaffian method~\parencite{Kasteleyn1963,Temperley1961,Fisher1961},
no solution has been found when unoccupied lattice sites (monomers) are allowed.
This stark dichotomy has captivated theoretical computer scientists
who focused on computational complexity theory,
and it was shown that 
the enumeration of the number of arrangements of dimers on quadratic lattices
belongs to \#P-complete class~\parencite{Jerrum1987}
and is considered intractable~\parencite{valiantComplexityEnumerationReliability1979,valiantComplexityComputingPermanent1979}.
For counting problems,
the class \#P plays the analogous role
of the more familiar NP class does for decision problems
(such as the well-known satisfiability problem).
The \#P-complete problems are at least as hard as the NP-complete
problems in computational complexity hierarchy.
Currently 
``P versus NP'' problem is perhaps the major
outstanding problem in theoretical computer science.

In this paper we demonstrate the structure of the number of arrangements
for a fixed number of dimers on two-dimensional quadratic lattices.
To our best knowledge this is the first time the structure of the solution to a \#P-complete problem
has been explicitly unveiled,
and the results show that computational complexity of calculating the number of covering $s$ dimers
on the lattice is $O(2^s)$.

\section{The nesting structure}

For a two-dimensional $n \times m$  quadratic lattices
with periodic boundary conditions in both directions,
we proved that canonical partition function
of $s$ rigid hard linear rods with length $k$
is a polynomial function of $N = n \times m$ with a degree of $s$,
provided that $m > (k-1)s + 1$ and $n > (k-1)s + 1$~\parencite{Kong2026}.

For monomer-dimer problem, denote the number of arrangements for $s$ dimers
on a lattice with $N$ sites as $a(N, s)$.  The above results then specializes
that $a(N, s)$ is a polynomial function of $N$
when $m > s+1$ and $n > s + 1$.
\begin{equation} \label{E:ac}
a(N, s) = c_0 N^s + c_1 N^{s-1} + \cdots + c_{s-1} N
= \sum_{i=0}^{s-1} c_i N^{s-i} .
\end{equation}

Using exact calculations
on two-dimensional lattices with cyclic boundary conditions in both directions,
we obtain $a(N, s)$ for $s$ up to $17$.
A partial results are show in Table~\ref{T:N} for $s=1, \dots, 9$.

\begin{table}[thbp] 
  \centering
  \caption{Polynomials in $N$ for $1 \le s \le 9$.}\label{T:N}
  {\tiny
    \begin{tabular}{RRRRRRRRRRRRRR}
      \hline\hline
      s  & N^{s}  & N^{s-1}  & N^{s-2}  & N^{s-3}  & N^{s-4}  & N^{s-5}  & N^{s-6}  & N^{s-7}  & N^{s-8} \\
      \hline
1 & 2 &    \\ 
2 & 2 & -7 &    \\ 
3 & 4/3 & -14 & 116/3 &    \\ 
4 & 2/3 & -14 & 611/6 & -521/2 &    \\ 
5 & 4/15 & -28/3 & 379/3 & -2375/3 & 9812/5 &    \\ 
6 & 4/45 & -14/3 & 905/9 & -2239/2 & 584627/90 & -47644/3 &    \\ 
7 & 8/315 & -28/15 & 526/9 & -1003 & 450641/45 & -833582/15 & 945688/7 &    \\ 
8 & 2/315 & -28/45 & 1199/45 & -5840/9 & 3511207/360 & -16379413/180 & 137529139/280 & -4768953/4 &    \\ 
9 & 4/2835 & -8/45 & 1346/135 & -4892/15 & 3688789/540 & -2827043/30 & 9529875731/11340 & -400990313/90 & 97398764/9 &    \\ 
      \hline
    \end{tabular}
    }
\end{table}

As a check of the correctness of the ingester calculations,
we compare our data with those derived from extensive data from Ising-model studies obtained using different methods~\parencite{Gaunt1969}.
In the thermodynamic limit 
the configurational grand-canonical partition $\Xi(z) = \lim_{N \to \infty} \Xi_N(z)^{1/N} = 1 + \sum_{s=1}^\infty  a(s) z^s $,
where it can be shown that $a(s) = a(N, s)|_{N=1}$.
Substitute $N=1$ in our $a(N, s)$ leads to Table~\ref{T:gaunt},
which agrees exactly with the square lattice column of Table~1 of~\textcite{Gaunt1969}
(Gaunt's table is up to $s=15$).

\begin{table}[thbp] 
  \centering
  \caption{
    The table for $a(1, s)$. See~\textcite{Gaunt1969}.
  }\label{T:gaunt}
  \begin{tabular}{RR}
    \hline\hline
 1 & 2 \\ 
 2 & -5 \\ 
 3 & 26 \\ 
 4 & -172 \\ 
 5 & 1288 \\ 
 6 & -10409 \\ 
 7 & 88594 \\ 
 8 & -782929 \\ 
 9 & 7119294 \\ 
 10 & -66205461 \\ 
 11 & 626921314 \\ 
 12 & -6025780853 \\ 
 13 & 58647507846 \\ 
 14 & -576914565407 \\ 
 15 & 5727408719806 \\ 
 16 & -57315706544813 \\ 
 17 & 577612294238894 \\
 \hline
\end{tabular}
\end{table}

From the table of $a(N, s)$ (the first few rows are presented in Table~\ref{T:N}),
there are clear and generalizable patterns for each column,
which contains the coefficients $c_i(s)$ of $N^{s-i}$ of $a(N, s)$.
These patterns of $c_i(s)$ can unambiguously established:
multiplied by a factor function of $s$ and $i$,
they themselves are polynomials
in $s$ with a degree of $i-1$ when $i>0$,
whose coefficients depend on $i$ but not $s$.

Denote the factor as $f_c(s, i)$. The factor $f_c(s, i)$ has the following expression
\begin{equation} \label{E:fc}
f_c(s, i) = \frac{(-1)^i 2^{s - 2i}}{(s - i - 1)!} .
\end{equation}
The first few columns for a given arbitrary $s$ are shown in Eq.~\eqref{E:c} (as coefficients for $N^s, N^{s-1}, \dots$)
We can reliably obtain $c_i(s)$ for $0 \le i \le 12$.

{\tiny
  \begin{equation} \label{E:c}
\resizebox{\textwidth}{!}{$
\begin{aligned}
  c_0   &= \frac{2^s}{s!} = \frac{2^s}{(s-1)!} \left[ \frac{7^0}{0!} s^{-1} \right], \\
  c_1 &= - \frac{2^{s-2}}{(s-2)!} \left[ \frac{7^1}{1!} s^0  \right], \\
  c_2 &= \frac{2^{s-4} } {(s-3)!} \left[ \frac{7^2}{2!} s + \frac{23}{6}\right] ,\\
  c_3 &= -\frac{2^{s-6} } {(s-4)!} \left[ \frac{7^3}{3!} s^2 + \frac{7}{1!} \cdot \frac{23}{6} s + 20\right] , \\
  c_4 &= \frac{2^{s-8} } {(s-5)!} \left[ \frac{7^4}{4!} s^3
    + \frac{7^2}{2!} \cdot \frac{23}{6} s^2
    + \frac{103^2}{3 \cdot 4!} s + \frac{19681}{180} \right] , \\
  c_5 &= -\frac{2^{s-10} } {(s-6)!} \left[ \frac{7^5}{5!} s^4
        + \frac{7^3}{3!} \cdot \frac{23}{6} s^3
        + \frac{7 \cdot 5569}{3 \cdot 1! \cdot 4!} s^2
        + \frac{13 \cdot 89 \cdot 131}{180} s
        + 708
        \right] ,\\
  c_6 &= \frac{2^{s-12} } {(s-7)!} \left[ \frac{7^6}{6!} s^5
        + \frac{7^4}{4!} \cdot \frac{23} {6} s^4
        + \frac{7^2 \cdot 3889} {3 \cdot 2! \cdot 4!} s^3
        + \frac{ 20897077 } {3^2 \cdot 6!} s^2
        + \frac{  79 \cdot 199 \cdot  383 } { 3! \cdot 180 } s
        + \frac{ 29 \cdot 43 \cdot 11689 } { 3^4 \cdot 5 \cdot 7 }
        \right], \\
  c_7 &= -\frac{2^{s-14} } {(s-8)!} \left[ \frac{7^7}{7!} s^6
        + \frac{7^5}{5!} \cdot \frac{23} {6} s^5
        + \frac{7^3 \cdot 3049} {3 \cdot 3! \cdot 4!} s^4
        + \frac{7 \cdot 29 \cdot 261581} {9 \cdot 6!} s^3
        + \frac{ 23573021 } {9 \cdot 5!} s^2
        + \frac{ 2 \cdot 8280499} {3^4 \cdot 5} s
        + \frac{ 2^3 \cdot 5 \cdot 3019} {3}
        \right] .
\end{aligned}
$}
\end{equation}
}

Denote the columns of $c_i(s)$ inside the square brackets in Eq.~\eqref{E:c}) by $d'_j(i)$ for $0 \le j \le i-1$
(the prime indicates that the coefficients differ from the unprimed ones by a factor $f_c(s, i)$),
then we have, for $i>0$,
    \begin{equation} \label{E:cd}
    c_i(s) = f_c(s, i)
    \left[
      d'_0 s^{i-1} + d'_1 s^{i-2} + \cdots + d'_{i-1}
      \right]
    = f_c(s, i)
    \sum_{j=0}^{i-1} d'_j(i) s^{i-1-j}
    = \sum_{j=0}^{i-1} d_j(s, i) s^{i-1-j}.
    \end{equation}

The columns of $c_i(s)$ (as shown in Eq.~\eqref{E:c})
again unveil clear patterns for the coefficients $d'_j(i)$:
they are polynomials in $i$ when $j>0$, with coefficients depending on $j$ but not $i$,
multiplied by a factor function of $i$ and $j$.

The first few columns of $c_i(s)$ inside brackets are listed in Eq.~\eqref{E:d}.
We can reliably obtain $d_j(s, i)$ for $0 \le j \le 7$.
Note due to the factor $f_c(s, i)$ outside the square bracket of Eq.~\eqref{E:c},
$d_j(s, i)$ is a function of $s$, $i$, and $j$,
but the coefficients inside the square bracket are independent of $s$:
\[
d_j(s, i) = f_c(s, i) d'_j(i).
\]
{\scriptsize
  \begin{equation} \label{E:d}
    \begin{aligned}
d'_0 =& \frac{7^i} {i!} %
    = \frac{7^i} {(i-1)!}  \left[ \frac{ 23^0  } {0! 6^0} i^{-1} \right] ,\\
d'_1 =& \frac{7^{i-2} \cdot 23 } {3! (i-2)!}
   = \frac{7^{i-2}  } {(i-2)!} \left[ \frac{ 23^1  } {1! 6^1} i^0 \right],  \\
d'_2 =& \frac{7^{i-4} } { (i-3)! } \left[ \frac{ 23^2 }{ 2! 6^2} i
      + \frac{1}{0!}\frac{ 2831 }{ 24}  \right], \\ %
d'_3  =& \frac{7^{i-6} } {  (i-4)! }
    \left[
    \frac{ 23^3 } { 3! 6^3} i^2
    + \frac{1}{1!} \frac{23}{6} \frac{ 2831}{24} i
    +  \frac{1376474 }{405}
    \right], \\    %
d'_4 =& \frac{7^{i-8} }  {  (i-5)! }
     \left[
    \frac{ 23^4 } { 4! 6^4} i^3
    + \frac{1}{2!} \frac{23^2}{6^2} \frac{ 2831 } { 24} i^2
    + \frac{ 3108135463 } { 155520 } i
    + \frac{ 622730773 } { 5184 }
           \right], \\ %
d'_5 =& \frac{7^{i-10} }  {  (i-6)! }
             \left[
             {\frac {23^5}{5! 6^5}}\,{i}^{4}
             +{\frac{1}{3!}  \frac{23^3}{6^3} \frac {2831}{24}}\,{i}^{3}
             +{\frac {48186163777}{933120}}\,{i}^{2}+{\frac {133962805199}{155520}}\,i+{\frac {4106113501}{810}}
             \right], \\ %
d'_6  =&
              \frac{ {7}^{i-12} } { (i-7)!}
              \left[ {\frac {23^6}{6! 6^6}}\,{i}^{5}
                +{\frac{1}{4!} \frac{23^4}{6^4} \frac {2831 }{24 }}\,{i}^{4}
                +{\frac {2788923407557}{33592320}}\,{i}^{3}
                +{\frac {3350449095911}{1244160}}\,{i}^{2}+{\frac {3306963804437521}{83980800}}\,i
                +{\frac {277418240388709}{1166400}}
                \right],  \\
\end{aligned}
    \end{equation}
}

Denote the factor outside the bracket of E.~\eqref{E:d} by $f_d(i, j)$
\begin{equation} \label{E:fd}
f_d(i, j) = \frac{7^{i - 2j}}{(i-j-1)!},
\end{equation}
and the coefficients of $i^{k-1}$ inside the bracket by $e'_k(j)$, for $0 \le k \le j-1$,
we have, for $j>0$,
\[
d'_j(i) = f_d(i, j) 
\left[
  e'_0 i^{j-1} + e'_1 i^{j-2} + \cdots + e'_{j-1}
  \right]
= f_d(i, j) 
\sum_{k=0}^{j-1} e'_k(j) i^{j-1-k} .
\]

The columns inside the square brackets of $d'_j(i)$ (as shown in Eq.~\eqref{E:d})
again unveil a pattern for the coefficients $e'_k(j)$:
they are polynomials in $j$ when $k>0$,
with coefficients depending on $k$, multiplied by a factor of function of $j$ and $k$.
Denote this factor by $f_e(j, k)$,
\begin{equation} \label{E:fe}
f_e(j, k) = \left(\frac{23}{6} \right)^{j-2k} \frac{1}{(j-k-1)!}.
\end{equation} 

For $e'_k(j)$ we can reliably obtain up to $k \le 4$, as shown in Eq.~\eqref{E:e}.
{\scriptsize
  \begin{equation} \label{E:e}
    \begin{aligned}
  e'_0 &= \frac{23^j}{ 6^j \cdot j!} = \frac{23^j}{ 6^j \cdot (j-1)!}  \left[ \frac{2831^0}{ 0! \cdot 24^0 } j^{-1} \right], \\
  e'_1 &= \frac{23^{j-2}}{ 6^{j-2} \cdot (j-2)! }
  \left[ \frac{2831^1}{1! \cdot 24^1 }  j^0 \right], \\
  e'_2 &=
  \frac{23^{j-4}}{ 6^{j-4} (j-3)!}
  \left[
    \frac{2831^2}{2! \cdot 24^2} j
    - \frac{1219727477} { 155520 } 
    \right],\\
  e'_3 &=
  \frac{23^{j-6}}{6^{j-6} (j-4)!}
  \left[
    \frac{2831^3}{3! \cdot 24^3} j^2
    -\frac{2831}{24}  \frac{1219727477}{ 155520 } j
    +\frac{75267510192}{24^3 \cdot 5}
    \right],\\
  e'_4 &=
  \frac{23^{j-8}}{6^{j-8} (j-5)!}
  \left[
    {\frac {2831^4}{4! \cdot 24^4}}\,{j}^{3}
    -{\frac{1}{2!} \frac {2831^2}{24^2} \frac{1219727477}{155520}}\,{j}^{2}
    +{\frac {7701215608818361849}{24^4 \cdot 145800}}\,j
    -{\frac {748530113262144671}{24^4 \cdot 14580}}
\right] .
  \end{aligned}
\end{equation}
}

Again the columns of the coefficients $f'_l(k)$ of $e'_k(j)$ show a clear pattern.
They are polynomials in $k$ when $l>0$, with coefficients depending on $j$, multiplied by a factor function of $k$ and $l$.
Denote this factor by $f_f(k, l)$ (not confused with $f'_l(k)$),
\begin{equation} \label{E:ff}
f_f(k, l) =  \left(\frac{2831}{24} \right)^{k-2l} \frac{1}{(k-l-1)!} .
\end{equation}

For $f'_l(k)$ we can reliably obtain up to $l \le 2$,
\begin{equation} \label{E:f}
  \begin{aligned}
    f'_0  &= \frac{2831^k}{ 24^k \cdot k! }
    = \left( \frac{2831}{24} \right)^k     \frac{1}{ (k-1)!}
    \left[ \left(-\frac{1219727477}{24^2 \cdot 270} \right)^0  k^{-1} \right], \\
    f'_1  &= \left( \frac{2831}{24} \right)^{k-2}  \frac{1}{ (k-2)!}         
         \left[ \left(-\frac{1219727477}{24^2 \cdot 270} \right)^1  k^{0} \right], \\
    f'_2  &=    \left( \frac{2831}{24} \right)^{k-4}  \frac{1}{ (k-3)!}
    \left[
      \frac{1}{2!}\left(-\frac{1219727477}{24^2 \cdot 270} \right)^2 k + \frac{583425045407739911}{24^4 \cdot 270 \cdot 180}
\right],
\end{aligned}
\end{equation}
and in general, for $k>0$,
\[
e'_k(j) = f_e(j, k) 
\sum_{l=0}^{k-1} f'_l(k) j^{k-1-l} .
\]

Compare $d'_j(i)$  with
$e'_k(j)$ and with $f'_l(k)$, we see that as we go with the sequence of indices
of the factors $f_c(s, i)$, $f_d(i, j)$, $f_e(j, k)$, $f_f(k, l)$,
\[
(s, i) \to
(i, j) \to
(j, k) \to
(k, l) \to
\cdots,
\]
the sequence of numbers inside the factor functions goes as

\[
2  \to
7 \to
23/6 \to
2831/24 \to
-1219727477/155520
\to
\cdots .
\]
This process continues, with the same pattern.

\section{The solution}

Combine the results together, we have the solution of
arranging $s$ dimers on a quadratic lattice with $N$ sites as
{\scriptsize
\begin{equation} \label{E:aNs}
  \begin{aligned} 
  a(N, s) &= c_0 N^s + \sum_{i=1}^{s-1} c_i(s) N^{s-i} \\
      &= c_0 N^s + \sum_{i=1}^{s-1} d_0 s^{i-1} N^{s-i} +
      \sum_{i=1}^{s-1}  \sum_{j=1}^{i-1} d_{j}(i) s^{i-j-1}  N^{s-i}  \\
      &= c_0 N^s
      + \sum_{i=1}^{s-1} d_0 s^{i-1} N^{s-i}
      + \sum_{i=1}^{s-1}  \sum_{j=1}^{i-1} e_0  i^{j-1}  s^{i-j-1}  N^{s-i} 
      + \sum_{i=1}^{s-1}
      \sum_{j=1}^{i-1}
      \sum_{k=1}^{j-1}
      e_{k}(j)
      i^{j-k-1}
      s^{i-j-1}
      N^{s-i} \\
  &=
  c_0 N^s
  + \sum_{i=1}^{s-1} d_0 s^{i-1} N^{s-i}
  + \sum_{i=1}^{s-1}  \sum_{j=1}^{i-1} e_0  i^{j-1}  s^{i-j-1}  N^{s-i}
  + \sum_{i=1}^{s-1}  \sum_{j=1}^{i-1} \sum_{k=1}^{j-1} f_0  j^{k-1} i^{j-k-1}  s^{i-j-1}  N^{s-i}
  +
  \cdots\\
  \end{aligned}
\end{equation}
}%
The solution has $s$ terms, with increasing number of nesting sums.
For a fixed $s$,
because the upper limit of each sum is one less than the previous one, 
there are $s-1$ sums in the last term.

Note that in the final expression of the solution in Eq.~\eqref{E:aNs},
only terms with zero index appeared.
As show in Eqs.~\eqref{E:c}, \eqref{E:d}, \eqref{E:e}, and \eqref{E:f},
the primed version of these terms are
\[
2^s/s!,
\quad
7^i/i!,
\quad
(23/6)^j/j!, \cdots .
\]

To facilitate the exposition, we map the indices from
$s$, $i$, $j$, etc. to $v_i, i=0, 1, \dots$,
with $s \to v_0$, $i\to v_1$, etc.
The factors $f_c(s, i)$, $f_d(i, j)$, etc.
map to $f_1(v_0, v_1)$, $f_2(v_1, v_2)$, etc.,
and denote the sequence of numbers appearing in these factors
\begin{equation} \label{E:g}
2, 7, 23/6, 2831/24, -1219727477/155520, \cdots
\end{equation}
by $g_i, i=1, 2, \dots$.
Let $z_i, i=1, 2, \dots$, denote the sequence of
$c_0(s)$, $d_0(s, i)$, $e_0(s, i, j), \dots$,
as shown in Eq.~\eqref{E:aNs}.

By this notation,
for $x>y\geq1$ define the following weights, using the constants $g_i$ in
\eqref{E:g}:
\begin{equation}\label{eq:weights}
 f_1(x,y)=\frac{(-1)^y2^{x-2y}}{(x-y-1)!},\qquad
 f_i(x,y)=\frac{g_i^{x-2y}}{(x-y-1)!}\quad(i\geq2).
\end{equation}
For a strictly decreasing chain $v_0>v_1>\cdots>v_r\geq1$ put
\begin{equation}\label{eq:z}
 z_{r+1}(v_0,\ldots,v_r)
 =\frac{g_{r+1}^{v_r}}{v_r!}
  \prod_{k=0}^{r-1}f_{k+1}(v_k,v_{k+1}).
\end{equation}

\begin{theorem}[Complete nested representation]\label{thm:full}
For every $s\geq1$, $m,n>s+1$, and $N=mn$, the exact matching count is
\begin{align}\label{eq:full}
 a(N,s)={}&\frac{2^s}{s!}N^s\\
 &+\sum_{r=1}^{s-1}
   \sum_{s=v_0>v_1>\cdots>v_r\geq1}
   N^{s-v_1}z_{r+1}(v_0,\ldots,v_r)
   \left(\prod_{k=0}^{r-2}
      v_k^{v_{k+1}-v_{k+2}-1}\right)
   v_{r-1}^{v_r-1}.\notag
\end{align}
The empty product for $r=1$ is one.  Every index-chain coefficient
structure in this formula holds for arbitrary $s$, rather than only
for the finitely computed examples.
\end{theorem}
We defer the proof of Theorem~\ref{thm:full} to a forthcoming paper.

\section{The sequence $g_i$}

The sequence $g_i$ plays a pivotal role in the solution.
If the sequence $g_i$ with $i$ from $1$ to $s$ is known,
we can calculate $f_i(v_i, v_{i+1})$ for $i=1$ to $s$,
which in turn determines $z_i(v_1, v_2, \dots, v_i)$, hence $a(N, s)$ by Eq.~\eqref{E:aNs}.
Conversely, if $a(N, s)$ are known,
$g_s$ can be obtained from Eq~\eqref{E:aNs}.
The first few expressions of $a(N, s)$ in terms of $g_i$
are shown in Eq.~\eqref{E_a_of_g}.

\begin{equation} \label{E_a_of_g}
  \begin{aligned}  
a(N, 1) &= g_{{1}}N  ,\\
a(N, 2) &= \frac{1}{2}\,{g_{{1}}}^{2}{N}^{2}-g_{{2}}N, \\
a(N, 3) &= \frac{1}{6}\,{g_{{1}}}^{3}{N}^{3}-g_{{2}}g_{{1}}{N}^{2}
+ \frac{1}{2}\,{\frac { \left( 3\,{g_{{2}}}^{2}+2\,g_{{3}} \right) N}{g_{{1}}}}, \\
a(N, 4) &= \frac{1}{24}\,{g_{{1}}}^{4}{N}^{4}
-\frac{1}{2}\,g_{{2}}{g_{{1}}}^{2}{N}^{3}+ \left( g_{{3}}+2\,{g_{{2}}}^{2} \right) {N}^{2}
-\frac{1}{6}\,{\frac { \left( 24\,g_{{3}}{g_{{2}}}^{2}+9\,{g_{{3}}}^{2}+16\,{g_{{2}}}^{4}+6\,g_{{4}} \right) N
  }{g_{{2}}{g_{{1}}}^{2}}}
\end{aligned}
\end{equation}

The list of $g_i$ for $1 \le 17$ is shown
in Table~\ref{T:g}.
The absolute values of these $g_i$ numbers grow
explosively.
A fit to the $g_i$ for $10 \le i \le 17$ as a function of $i$
shows that sequence grows double-exponentially:
\begin{equation} \label{E:loglog}
\log( \log( |g_i| ) ) \approx 0.693125   i -1.28179 .
\end{equation}

\renewcommand{\arraystretch}{2.2} %
\begin{table}[thbp] 
  \centering
  \caption{
    The sequence of $g_i$ for $1 \le i \le 17$.
  }\label{T:g}
  \begin{tabular}{RRRR}
    \hline\hline
    i & \text{ exact } g_i  & g_i \text{ in float } & \log(\log(|g_i|)) \\
    \hline
1 & 	2 & 	 \quad  \num{  2. } & 	 \quad  \num{  -.366513 }  \\ 
2 & 	7 & 	 \quad  \num{  7. } & 	 \quad  \num{  .665730 }  \\ 
3 & 	\dfrac{23}{6} & 	 \quad  \num{  3.83333 } & 	 \quad  \num{  .295449 }  \\ 
4 & 	\dfrac{2831}{24} & 	 \quad  \num{  117.958 } & 	 \quad  \num{  1.56242 }  \\ 
5 & 	-\dfrac{1219727477}{155520} & 	 \quad  \num{  -7842.90 } & 	 \quad  \num{  2.19359 }  \\ 
6 & 	\dfrac{583425045407739911}{16124313600} & 	 \quad  \num{  .361829e8 } & 	 \quad  \num{  2.85671 }  \\ 
7 & 	-\dfrac{571[\dots32 \text{ digits}\dots]961}{140[\dots17 \text{ digits}\dots]000} & 	 \quad  \num{  -.406950e16 } & 	 \quad  \num{  3.58191 }  \\ 
8 & 	-\dfrac{577[\dots69 \text{ digits}\dots]961}{131[\dots39 \text{ digits}\dots]000} & 	 \quad  \num{  -.439426e31 } & 	 \quad  \num{  4.25643 }  \\ 
9 & 	\dfrac{273[\dots147 \text{ digits}\dots]973}{326[\dots86 \text{ digits}\dots]000} & 	 \quad  \num{  .838551e61 } & 	 \quad  \num{  4.94365 }  \\ 
10 & 	-\dfrac{245[\dots301 \text{ digits}\dots]049}{710[\dots177 \text{ digits}\dots]000} & 	 \quad  \num{  -.345719e124 } & 	 \quad  \num{  5.65059 }  \\ 
11 & 	-\dfrac{561[\dots610 \text{ digits}\dots]491}{299[\dots363 \text{ digits}\dots]000} & 	 \quad  \num{  -.187540e248 } & 	 \quad  \num{  6.34453 }  \\ 
12 & 	-\dfrac{296[\dots1224 \text{ digits}\dots]161}{598[\dots731 \text{ digits}\dots]000} & 	 \quad  \num{  -.495049e493 } & 	 \quad  \num{  7.03392 }  \\ 
13 & 	-\dfrac{107[\dots2458 \text{ digits}\dots]901}{965[\dots1470 \text{ digits}\dots]000} & 	 \quad  \num{  -.110885e988 } & 	 \quad  \num{  7.72875 }  \\ 
14 & 	-\dfrac{332[\dots4919 \text{ digits}\dots]969}{126[\dots2945 \text{ digits}\dots]000} & 	 \quad  \num{  -.262360e1975 } & 	 \quad  \num{  8.42206 }  \\ 
15 & 	-\dfrac{153[\dots9849 \text{ digits}\dots]287}{146[\dots5901 \text{ digits}\dots]000} & 	 \quad  \num{  -.104828e3949 } & 	 \quad  \num{  9.11500 }  \\ 
16 & 	\dfrac{251[\dots19703 \text{ digits}\dots]311}{142[\dots11807 \text{ digits}\dots]000} & 	 \quad  \num{  .176359e7897 } & 	 \quad  \num{  9.80818 }  \\ 
17 & 	-\dfrac{157[\dots39415 \text{ digits}\dots]117}{930[\dots23622 \text{ digits}\dots]000} & 	 \quad  \num{  -.168743e15793 } & 	 \quad  \num{  10.5013 }  \\ 
\hline
  \end{tabular}
\end{table}
\renewcommand{\arraystretch}{1.0} %

The $g_i$ sequence has both positive and negative terms.
While double-exponential sequences occasionally surface in highly contrived mathematical settings,
they are exceptionally rare in natural problems.
Within enumerative combinatorics and statistical physics,
such behavior is virtually unknown.
To the best of our knowledge, the sequence discovered here
represents the first reported occurrence of its kind in this domain.

The number of digits in the numerators and denominators of $g_i$ grows rapidly,
much faster than the coefficients of $N^i$ in $a(N,s)$.
For $a(N, s)$, $g_j$ only appears in the coefficients of $N^i$ when $i \le s-j+1$.
Compare with the coefficients of $a(N, s)$, $g_i$ has much bigger number of digits.
For example, when $s=17$, the coefficient of $N$ in $a(N, 17)$ is $14719835348283940/17$,
but the numerator and denominator of $g_{17}$ have $39421$ and $23628$ digits,
respectively (see Table~\ref{T:g}).

The advantage of expressing $a(N,s)$ in the sequence of $g_i$ is that
the elements in the sequence are independent of $s$;
the same sequence works for arbitrary $s$. 
Once $g_i$ is known up to $s$, for bigger $s$ such as $s+1$,
there is only one unknown element $g_{s+1}$ to be found, and this can be found 
computationally for any $N > (s+2)^2$. 
Actually the equations shown in Eq~\eqref{E_a_of_g} hold
formally
for any $N > 0$.
The Ising-related method can be used here to get the value of $a(N, s)$ for $N=1$,
then the values of $g_{s+1}$ can be obtained.

\section{Computational complexity}

Each term in the expression Eq.~\eqref{E:aNs} represents a specific loop depth
$t$. The number of operations required to evaluate a term with $t$ nested loops
is given by the binomial coefficient $\binom{s-1}{t}$.
If we add the operations for all terms together, we get the total operations
\[
\sum_{t=0}^{s-1} \binom{s-1}{t} = 2^{s-1},
\]
hence
the computational complexity is $O(2^s)$.

The first term of Eq.~\eqref{E:aNs} is
\[
\frac{2^s}{s!} N^s ,
\]
and
the second term of the sum of $d_0(s, i)$ has a closed form
\[
{\frac {N \left( 4\, \left( 4\,N-7\,s \right) ^{s-1}-{N}^{s-1}{4}^{s}
 \right) }{s!\,{2}^{s}}} ,
\]
which gives the $N^{s-1}$ term in $a(N,s)$ in addition to other terms.

The sum of $e_0(s, i, j)$ contains a double sums, which can be simplified into a single sum

\[
\sum _{i=1}^{s-1}
\frac{7\, \left( -1 \right) ^{i}{N}^{s-i}{2}^{s-2\,i} }{(s-1)!}
 \left( 7\,s+{\frac {23}{42}}\,i \right) ^{i-1}{s-1\choose i}  .
\]
This expression is not summable.
For $e_0(s, i, j)$ term and up, the summands contains $i^j$, $j^k$, $i^k$, etc.
These summands are not hypergeometric functions and are not summable to
simpler forms.

Given the uniqueness of power series expansions, the solution shown in Eq.~\eqref{E:aNs} is in its simplest form;
consequently, our results demonstrate that the computational complexity of this \#P-complete problem is exponential.

\newgeometry{margin=1cm} %
\begin{landscape}
\thispagestyle{empty}
\newgeometry{top=0.1in, bottom=0.1in, left=0.1in, right=0.1in}

\restoregeometry
\end{landscape}

\printbibliography
\end{document}